# A little noticed right triangle

*Konstantine 'Hermes' Zelator*
*Department of Mathematics*
*College of Arts and Sciences*
*Mail Stop 942*
*University of Toledo*
*Toledo, OH 43606-3390*
*U.S.A.*

## 1. Introduction

Given a right triangle $AB\Gamma$, there is another right triangle which naturally arises from $AB\Gamma$, and which seems to have escaped much attention. In **Figure 1**, a *nonisosceles* triangle (with right angle at the vertex *A*) is depicted. Let *O* be the midpoint of the hypotenuse $\overline{B\Gamma}$, which of course is also the center of the triangles' circumscribed circle. Consider the line segment $\overline{AO}$; it divides the right triangle $AB\Gamma$ into two noncongruent isosceles triangles $ABO$ and $A\Gamma O$, each having area $\frac{1}{2}E$, where E is the area of triangle $AB\Gamma$. Let $O_1$ and $O_2$ be respectively the centers of the circumscribed circles of triangles $ABO$ and $A\Gamma O$.

It is triangle $OO_1O_2$, which is the subject matter of this work. Note that since $AB\Gamma$ is a nonisosceles triangle, one of the points $O_1$ and $O_2$ must lie in the exterior of $AB\Gamma$, while the other in its interior. Specifically if length $\left|\overline{A\Gamma}\right| = \beta > \left|\overline{AB}\right| = \gamma$, then $O_2$ lies in the exterior of $AB\Gamma$, while $O_1$ in its interior (as in Figure 1 – while when $\beta < \gamma$, it is the other way around).

In Section 3, we easily show by using a few simple geometric arguments, that $OO_1O_2$ is indeed a right triangle, similar to triangle $AB\Gamma$.

In Section 4, we computationally establish formulas for the lengths $\left|\overline{OO_1}\right|, \left|\overline{OO_2}\right|$, and $\left|\overline{O_1O_2}\right|$; in terms of triangle $AB\Gamma$'s sidelengths $\alpha, \beta, \gamma$. These are simple rational expressions. We also compute the area of $OO_1O_2$ in terms of $\alpha, \beta, \gamma$.

In Section 2, we derive the formula $R = \frac{\alpha\beta\gamma}{4E}$, for any triangle $AB\Gamma$; where $R$ is the radius of the circumscribed circle, and $E$ the triangle's area (see **Figure 2**).

We make use of this result in Section 4.

In Section 5, we take a look at the nearby trapezoid $O_1M_1M_2O_2$, where $M_1$ is the midpoint of the segment $\overline{OB}$, and $M_2$ the midpoint of $\overline{O\Gamma}$. We calculate the trapezoid's four sidelengths, as well as its two diagonal lengths $d_2 = \left|\overline{O_2M_1}\right|$ and $d_1 = \left|\overline{O_1M_2}\right|$; and also this trapezoid's area.

In Section 6, some number theory enters the picture. We consider the case when $AB\Gamma$ is a Pythagorean triangle, and we examine both the triangle $OO_1O_2$ and the trapezoid in that case. Finally, in Section 7 we offer three numerical examples.

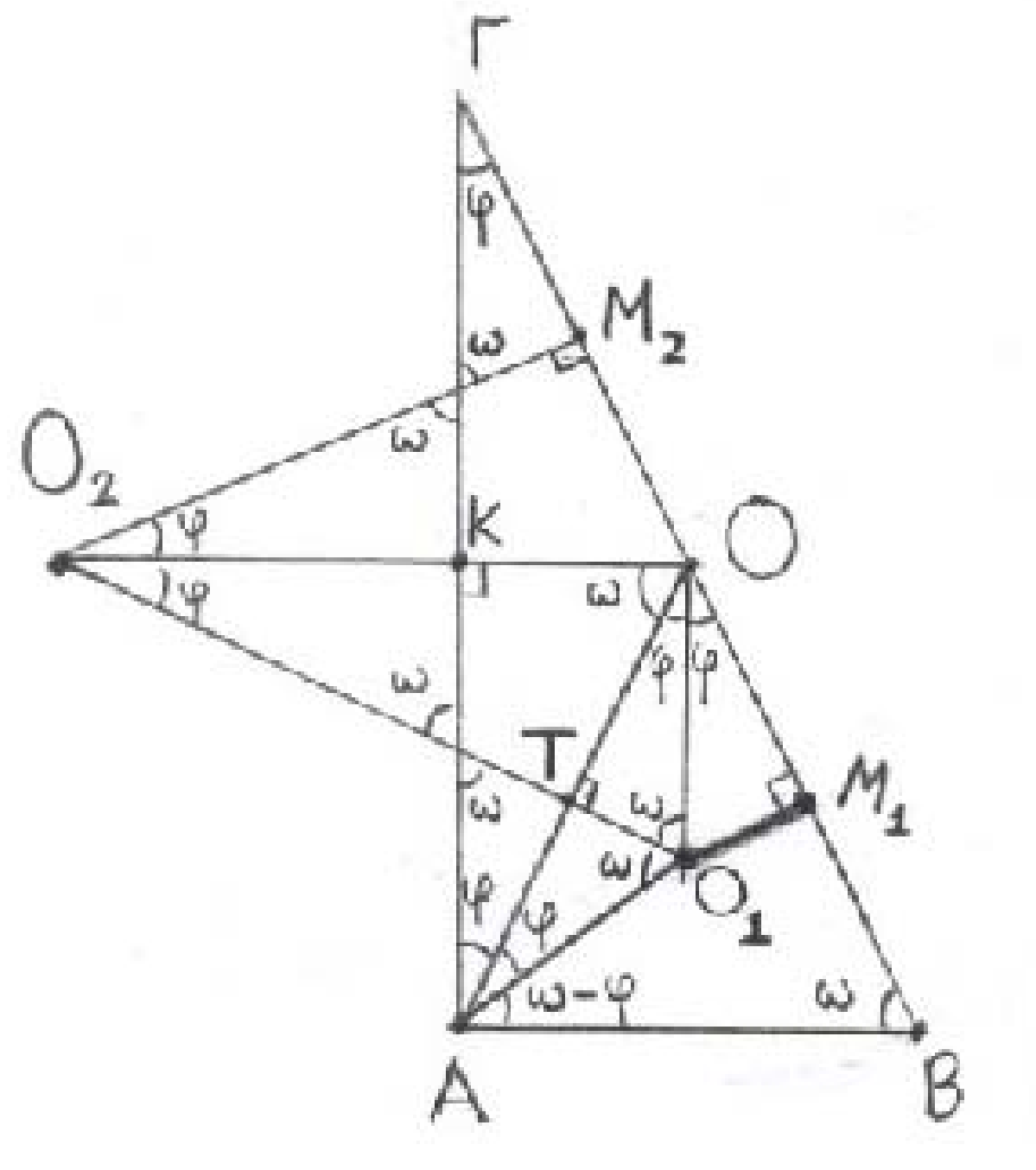

**Figure 1**

**2. A derivation of the formula** $R = \frac{\alpha\beta\gamma}{4E}$

In **Figure 2**, a triangle $AB\Gamma$ is inscribed in a circle of center *O* and radius *R*. If $\theta$ is the degree (or radian) measure of the triangle's internal angle at the vertex *A*; and also, with $|\overline{B\Gamma}| = \alpha$, $|\overline{BB_1}| = 2R$ ($B_1$ is the point on the circle, antidiametrical of $B$), $R = |\overline{O\Gamma}| = |\overline{OA}| = |\overline{OB}|$; we have from the right triangle $B\Gamma B_1$, $2R\sin\theta = \alpha$.

Moreover, let $B_2$ be the foot of the perpendicular drawn from *B* to the side $\overline{A\Gamma}$; and $h_\beta = |\overline{BB_2}|$. Then, $h_\beta = \gamma \sin\theta$; $\sin\theta = \frac{h_\beta}{\gamma}$. And so we obtain,

$$\alpha = 2R\sin\theta = 2R \cdot \frac{h_\beta}{\gamma} \quad \Leftrightarrow \quad 2Rh_\beta = \alpha \cdot \gamma .$$

Also, $E = \frac{1}{2}\beta \cdot h_\beta$, where $\gamma = |\overline{AB}|$, $\beta = |\overline{A\Gamma}|$, and E stands for the area of the triangle $AB\Gamma$. Solving for $h_\beta$ in terms of E gives $h_\beta = \frac{2E}{\beta}$; and thus $2Rh_\beta = \alpha\gamma$ further yields

$$2R \cdot \left(\frac{2E}{\beta}\right) = \alpha\gamma \quad \Leftrightarrow \quad \boxed{R = \frac{\alpha\beta\gamma}{4E}} \qquad (1)$$

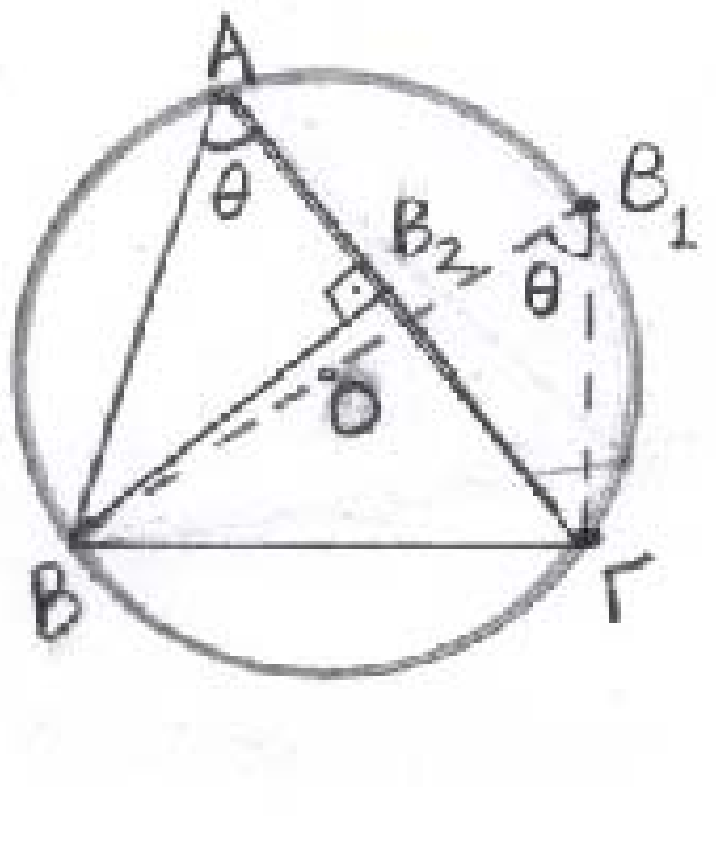

**Figure 2**

**3. The right triangle $OO_1O_2$**

Now we go back to **Figure 1**. It is easy to see that $O_1OO_2$ is a right triangle, in fact similar to $AB\Gamma$. Since $O_1$ is the center of the circumscribed circle of the isosceles triangle $AOB$, with $|\overline{AO}| = |\overline{OB}|$, it follows that $O_1$ lies on the bisector of the angle $< AOB$, since that angle bisector coincides with the perpendicular bisector of $\overline{AB}$. In other words, $\overline{OO_1}$ is perpendicular to $\overline{AB}$. Thus, if $\varphi$ is the degree measure of the angle $< A\Gamma B$, then since $\overline{OO_1}$ is parallel to $\overline{\Gamma A}$ (both being perpendicular to $\overline{AB}$), the angles $< A\Gamma B$ and $< O_1OB$ must have the same degree measure $\varphi$, based on the relevant theorem of Euclidean geometry which postulates that if two parallel straight lines are cut by another transversal line, then any pair of two angles, one internal the other external, and lying on the same side of the transversal; must have the same degree measure. Now, since $< O_1OB$ has degree measure $\varphi$, so does the angle $< AOO_1$, since as we have explained above, the segment $\overline{OO_1}$ lies on the bisector of the angle $< AOB$. Moreover, since the segment $\overline{OO_2}$ lies, by definition, on the perpendicular bisector of $\overline{A\Gamma}$; it follows that $\overline{OO_2}$ is parallel to $\overline{AB}$ (since $\overline{AB}$ is perpendicular to $\overline{A\Gamma}$ as well). Applying again a fundamental theorem of Euclidean geometry, which postulates that if two straight parallel lines are cut by a transversal, then any pair of interior alternate angles (i.e. lying on opposite sides of the transversal), must have the same degree measure. Let $\omega$ be the degree measure of angle $< OAB$. Since the degree measure of the angle $< OAB$ is $\omega$ (same as the measure of the angle $< OBA$), it then follows that the degree measure of the angle $< AOO_2$ is also $\omega$. Finally, since the line segment $\overline{OA}$ lies in the interior of the angle $< O_2OO_1$, it follows that the degree measure of the angle $< O_2OO_1$ is equal to the sum of degree measures of the angles $< O_2OA$ and $< AOO_1$; hence it is equal to $\omega + \varphi = 90^\circ$, which proves that $O_2OO_1$ is a right triangle with right angle at $O$. The triangle $OO_1O_2$ is similar to the triangle $AB\Gamma$, since it is easily seen that the angle $< O_2O_1O$ has measure $\omega$, while the angle $< OO_2O_1$ has measure $\varphi$.

*Immediate Observations*

With $\left|\overline{\Gamma B}\right| = \alpha$, $\left|\overline{A\Gamma}\right| = \beta$, $\left|\overline{AB}\right| = \gamma$; we have

$$\left|\overline{\Gamma M_2}\right| = \left|\overline{M_2 O}\right| = \left|\overline{OM_1}\right| = \left|\overline{M_1 B}\right| = \alpha/4 = \left|\overline{AT}\right| = \left|\overline{TO}\right|$$

$$\left|\overline{AO}\right| = \left|\overline{\Gamma O}\right| = \left|\overline{BO}\right| = \alpha/2$$

$$\left|\overline{\Gamma K}\right| = \left|\overline{KA}\right| = \beta/2 .$$

The following seven angles have the same degree measure $\varphi$:

$$< A\Gamma B,\ < AOO_1,\ < O_1OB,\ < OAO_1,\ < \Gamma AO,\ < M_2O_2O, \text{ and } < O_1O_2O .$$

On the other hand the following angles have degree measure $\omega$:

$$< \Gamma BA,\ < OAB,\ < AO_1O_2,\ < O_2O_1O,\ < O_2OA .$$

Also, we have $R_1 = \left|\overline{O_1O}\right| = \left|\overline{O_1A}\right| = \left|\overline{O_1B}\right|$ and $R_2 = \left|\overline{O_2O}\right| = \left|\overline{O_2\Gamma}\right| = \left|\overline{O_2A}\right|$, where $R_1$ is the radius of triangle $AOB$'s circumscribed circle; $R_2$ the radius of triangle $AO\Gamma$'s circumscribed circle.

**4. Computing the lengths $R_1, R_2$, and $\left|\overline{O_1O_2}\right|$**

If $E_1$ and $E_2$ are the areas of the isosceles triangles $AOB$ and $AO\mathrm{T}$ respectively; and $E$ the area of the triangle $AB\Gamma$ then,

$$\boxed{E_1 = E_2 = \frac{1}{2}E = \frac{\beta\gamma}{4}} \quad (2)$$

Combining (1) and (2) we see that,

$$\left|\overline{OO_1}\right| = R_1 = \frac{\left|\overline{OA}\right|\left|\overline{OB}\right|\left|\overline{AB}\right|}{4E_1} = \frac{\left(\alpha/2\right)\cdot\left(\alpha/2\right)\cdot\gamma}{\beta\gamma};$$

similarly,

$$\left|\overline{OO_2}\right| = R_2 = \frac{\left|\overline{OA}\right|\left|\overline{O\Gamma}\right|\left|\overline{A\Gamma}\right|}{4E_2} = \frac{\left(\alpha/2\right)\cdot\left(\alpha/2\right)\cdot\beta}{\beta\gamma}.$$

Altogether,

$$\boxed{\begin{aligned} R_1 &= \frac{\alpha^2}{4\beta} \\ R_2 &= \frac{\alpha^2}{4\gamma} \end{aligned}} \quad (3)$$

Let $x = \left|\overline{O_1M_1}\right|$. From the congruent right triangles $O_1TO$ and $O_1M_1O$, we see that $\left|\overline{O_1M_1}\right| = \left|\overline{O_1T}\right|$; and so $\left|\overline{O_1T}\right| = x$ as well. Similarly we have $y = \left|\overline{O_2M_2}\right| = \left|\overline{O_2T}\right|$.

Consequently,

$$\left|\overline{O_1O_2}\right| = \left|\overline{O_1T}\right| + \left|\overline{TO_2}\right| = x + y \qquad (4)$$

From the right triangle $M_1O_1O$, $x = \sqrt{R_1^2 - \left(\alpha/4\right)^2}$ by (3) $= \sqrt{\frac{\alpha^4}{16\beta^2} - \frac{\alpha^2}{16}} = \frac{\alpha}{4\beta}\cdot\sqrt{\alpha^2 - \beta^2}$ ;

$x = \frac{\alpha}{4\beta}\sqrt{\gamma^2}$ ; $x = \frac{\alpha\gamma}{4\beta}$ (since $\alpha^2 = \beta^2 + \gamma^2$). Working similarly with the triangle $M_2O_2O$ we find that $y = \frac{\alpha\beta}{4\gamma}$. Putting these two together,

$$\boxed{x = \frac{\alpha\gamma}{4\beta} \quad , \quad y = \frac{\alpha\beta}{4\gamma}} \quad (5)$$

From (4) and (5) $\Rightarrow$ $\left|\overline{O_1O_2}\right| = \frac{\alpha\gamma}{4\beta} + \frac{\alpha\beta}{4\gamma}$; which when combined with $\alpha^2 = \beta^2 + \gamma^2$

gives,

$$\boxed{\left|\overline{O_1O_2}\right| = \frac{\alpha^3}{4\beta\gamma}} \quad (6)$$

Also, the area of the right triangle $OO_1O_2$ is equal to $\frac{1}{2} \cdot R_1R_2$; which when combined with (3) yields,

$$\boxed{\text{Area of triangle } OO_1O_2 = \frac{\alpha^4}{32\beta\gamma}} \quad (7)$$

**5. The nearby trapezoid $O_1O_2M_2M_1$**

We have already computed, in terms of $\alpha, \beta$, and $\gamma$; three of the four sidelengths of this trapezoid. Namely, the lengths $\left|\overline{O_1O_2}\right|$, $\left|\overline{O_1M_1}\right| = x$, and $\left|\overline{O_2M_2}\right| = y$ according to formulas (5) and (6). The fourth sidelength is $\left|\overline{OM_1}\right| = \alpha/4$.

Moreover, the area of the trapezoid $O_1O_2M_2M_1$ is given by, area of $O_1O_2M_2M_1 = \frac{1}{2}(x+y) \cdot \left|\overline{M_1M_2}\right| = \frac{(x+y)\alpha}{4}$, and by using (5) and $\alpha^2 = \beta^2 + \gamma^2$ we arrive at,

$$\boxed{\text{Area of trapezoid } O_1O_2M_2M_1 = \frac{\alpha^4}{16\beta\gamma}} \quad (8)$$

There are two diagonals $\overline{M_1O_2}$ and $\overline{O_1M_2}$ of the trapezoid $O_1O_2M_2M_1$. Let $d_2$ and $d_1$ be their lengths respectively. From the right triangles $O_1M_2M_1$ and $O_2M_1M_2$ we obtain (since $\left|\overline{M_1M_2}\right| = \alpha/2$) $d_1 = \sqrt{x^2 + (\alpha/2)^2}$ and $d_2 = \sqrt{y^2 + (\alpha/2)^2}$ .

Implementing (5) further gives,

$$\boxed{d_1 = \frac{\alpha}{4\beta}\sqrt{\gamma^2 + 4\beta^2} \quad , \quad d_2 = \frac{\alpha}{4\gamma}\sqrt{\beta^2 + 4\gamma^2}} \qquad (9)$$

## 6. When $AB\Gamma$ is a Pythagorean triangle

In this section we venture a bit into number theory. If all three sidelengths $\alpha, \beta, \gamma$ are integers, then the triangle $AB\Gamma$ is what is commonly referred to as a Pythagorean triangle. Clearly then, formulas (3) and (6) tell us that the sidelengths $R_1$, $R_2$, $\left|\overline{O_1O_2}\right|$ of the triangle, while not necessarily integers, are certainly rational numbers.

Since $\{\alpha, \beta, \gamma\}$ is a Pythagorean triple, as it is well known (may refer to [1], [2], or [3]) we must have,

$$\alpha = \delta\left(m^2 + n^2\right) \quad , \quad \beta = \delta(2mn) \quad , \quad \gamma = \delta\left(m^2 - n^2\right) \quad , \qquad (10)$$

for some positive integers $\delta, m, n$ such that $(m, n) = 1$ (i.e. *m* and *n* are relatively prime), $m > n$, and $m + n \equiv 1$ (mod 2) (i.e. *m* and *n* have different parities – one of them is odd, the other even).

(Note: Obviously, $\beta = \delta(m^2 - n^2)$, $\gamma = \delta(2mn)$; is the other possibility. But there is no need to distinguish between two cases, since $\beta$ and $\gamma$ are interchangeable, there is no additional information about $\beta$ or $\gamma$ given here).

That the above parametric formulas (involving three parameters: $\delta$, $m$, and $n$) generate the entire family of Pythagorean triples is a well-known fact in number theory. Almost every undergraduate text or book in number theory, has some material on Pythagorean triangles (typically the parametric formulas and their derivation and some exercises on Pythagorean triples.)

Let us apply formulas (10) on (3), (6), and (7). After some algebra we obtain,

$$\boxed{\begin{array}{l} R_1 = \dfrac{\delta(m^2+n^2)^2}{8mn},\ R_2 = \dfrac{\delta(m^2+n^2)}{4(m^2-n^2)},\ \left|\overline{O_1O_2}\right| = \dfrac{\delta(m^2+n^2)^3}{8mn(m^2-n^2)} \\ \text{Area of triangle } OO_1O_2 = \dfrac{\delta^2(m^2+n^2)^4}{64mn(m^2-n^2)} \end{array}} \quad (11)$$

Likewise (5), (8), and (9) when combined with (10); they produce

$$\boxed{\begin{array}{l} x = \dfrac{\delta(m^2+n^2)(m^2-n^2)}{8mn},\ y = \dfrac{\delta mn(m^2+n^2)}{2(m^2-n^2)} \\ \text{Area of trapezoid } O_1O_2M_2M_1 = \dfrac{\delta^2(m^2+n^2)^4}{32mn(m^2-n^2)} \\ \text{And } d_1 = \dfrac{\delta(m^2+n^2)\sqrt{m^4+14m^2n^2+n^4}}{8mn},\ d_2 = \dfrac{\delta(m^2+n^2)\sqrt{m^4-m^2n^2+n^4}}{2(m^2-n^2)} \end{array}} \quad (12)$$

Now, let us consider the question: for what choice of the parameter $\delta$, are all three rationals $R_1$, $R_2$, $\left|\overline{O_1O_2}\right|$ actually integers? First note that since $(m,n)=1$ and $m+n\equiv 1(\text{mod}2)$, it follows that,

$$\left(\left(m^2+n^2\right)^{t_1}\ ,\ 8mn\left(m^2-n^2\right)^{t_2}\right)=1,$$

for any values of the nonnegative integer exponents $t_1$ and $t_2$. (This can be assigned as an exercise in an elementary theory course).

This then proves that,

$$\left(\left(m^2+n\right)^2, 8mn\right)=1=\left(\left(m^2+n^2\right)^2, 4\left(m^2-n^2\right)\right)=\left(\left(m^2+n^2\right)^3, 8mn\left(m^2-n^2\right)\right) \qquad (13)$$

Further, we also know from number theory that if an integer *a* divides the product $bc$ of integers *b* and *c*; and $(a,b)=1$. Then *a* must be a divisor of *c*. For this, refer to [1] or [3].

Which shows that in view of (13) and formulas (11), and in order for all three $R_1$, $R_2$, $\left|\overline{O_1O_2}\right|$ to be integers, it is necessary and sufficient that $\delta$ be divisible by all three integers $8mn$, $4\left(m^2-n^2\right)$, and $8mn\left(m^2-n^2\right)$. We need one more basic result from number theory.

An integer $\delta$ is divisible by the positive integers $a_1, a_2, \ldots, a_r$ if, and only if, $\delta$ is divisible by *L*, the least common multiple of $a_1, a_2, \ldots, a_r$ (*L* being the smallest positive integer divisible by the *r* numbers).

Applying the above with $a_1 = 8mn$, $a_2 = 4(m^2 - n^2)$, and $a_3 = 8mn(m^2 - n^2)$. Their least common multiple is $L = 8mn(m^2 - n^2)$ (the reader should verify).

The right triangle $OO_1O_2$ is a Pythagorean one, precisely when $\delta = K \cdot 8mn(m^2 - n^2)$, where $K$ can be any positive integer.

A Pythagorean triangle is called **primitive**, if its three sidelengths have no factor in common, other than 1. When $\delta = 1$, and only then, the triangle $AB\Gamma$ is primitive. In that case, all three sidelengths $R_1$, $R_2$, $|\overline{O_1O_2}|$ of triangle $OO_1O_2$ are proper rational numbers, none of them being an integer. On the other hand, if $\delta = K \cdot 8mn(m^2 - n^2)$, both right triangles $AB\Gamma$ and $OO_1O_2$ are nonprimitive; $OO_1O_2$ is nonprimitive since $(m^2 + n^2)^2$ is a common factor in the case of $R_1$, $R_2$, and $|\overline{O_1O_2}|$. In conclusion, if $OO_1O_2$ is a Pythagorean triangle, it is always nonprimitive.

Next, we apply similar considerations to the trapezoid $O_1O_2M_2M_1$, whose sidelengths are $x$, $y$, $|\overline{O_1O_2}|$, and $|\overline{M_1M_2}| = \frac{\alpha}{2} = \frac{\delta(m^2+n^2)}{2}$.

Using the same type of reasoning (as in the case of the triangle $OO_1O_2$), and by virtue of formulas (12) and conditions (13), we see once again that,

As in the case of triangle $OO_1O_2$, the trapezoid $O_1O_2M_2M_1$ will have all four of its sidelengths integers if, and only if, $\delta = K \cdot 8mn(m^2 - n^2)$, where K is a positive integer. Also, in this case, the trapezoid's area will also be an integer, as (12) clearly shows.

And this brings us to the last question in this paper.

What about the diagonal lengths $d_1$ and $d_2$ in (12), for $\delta = K \cdot 8mn(m^2 - n^2)$? Can either of them be an integer? The answer is no, neither $d_1$ nor $d_2$ can be integers when $\delta = K \cdot 8mn(m^2 - n^2)$. First, as we can easily see from (12), $d_1$ will be an integer for those integer values of $\delta$, if and only if, the real number $\sqrt{m^4 + 14m^2n^2 + n^4}$ is an integer. Similarly $d_2$ will be an integer precisely when $\sqrt{m^4 - m^2n^2 + n^4}$ is an integer.

Consider the following result from number theory.

| |
|---|
| Let *r* and $\ell$ be positive integers. Then $\sqrt[r]{\ell}$ (the *r*th root of $\ell$) is a rational number if, and only if, $\ell$ is the *r*th power of some positive integer *i*; that is $\ell = i^r$. Equivalently, $\sqrt[r]{\ell}$ is either an integer or an irrational number; the former occurring precisely when $\ell = i^r$. In particular, the square root (the case $r = 2$) of $\ell$ will be either an integer or an irrational number, the former occurring exactly when $\ell = i^2$, for some positive integer *i*. The interested reader may refer to [1] or [3]. |

Applying this to $d_1$ above, we see that $d_1$ will be integer precisely when $m^4 + 14m^2n^2 + n^4 = i^2$, for some positive integer *i*. Which means that $(m, n, i)$ would be a solution in positive integers, to the Diophantine equation $x^4 + 14x^2y^2 + y^4 = z^2$. According to L.E. Dickson's monumental book (see [2]), Euler was the first historically known individual to prove that all the solutions to the above equations are given by,

$$x = y = \delta \quad , \quad z = 4\delta^2 ,$$

where $\delta$ can be any positive integer.

Thus, under the condition $(x, y) = 1$, the above equation has only one solution, namely $x = y = 1, z = 4$. The integers *m* and *n* satisfy the same condition, $(m, n) = 1$. But *m* and *n* have different parities (one is even; the other odd); thus, $(m, n) \neq (1,1)$, which clearly shows that $(m, n, i)$ cannot really be a solution to the Diophantine equation $x^4 + 14x^2y^2 + y^4 = z^2$, and so, $\sqrt{m^4 + 14m^2n^2 + n^4}$ must be an irrational number.

We argue similarly in the case of $d_2$: in order for $d_2$ to be rational, it is necessary and sufficient that $\sqrt{m^4 - m^2n^2 + n^4} = i^2$, for some positive integer *i*; $m^4 - m^2n^2 + n^4 = i^2$, which means that $(m, n, i)$ would be a solution, in positive integers, to the Diophantine equation $x^4 - x^2y^2 + y^4 = z^2$.

L. Pocklington has proved (see [4]) , that all the solutions in positive integers, of the last equation, are given by $x = y = \delta, z = \delta^2$, where $\delta$ can be any positive integer. Under the condition $(x, y) = 1$, there is only one solution, $x = y = z = 1$. Since $(m, n) = 1$, and $m + n \equiv 1 (\text{mod } 2)$, once more we see that $(m, n, i)$ above cannot be a solution. Thus, $\sqrt{m^4 - m^2n^2 + n^4}$ must be an irrational number. Pocklington's proof is also succinctly presented in W. Sierpinski's book (see [3]).

When $AB\Gamma$ is a Pythagorean triangle, the two diagonal lengths $d_1$ and $d_2$ are always irrational numbers.

In a paper by S. Sastry (see [5]), a family of *Heron Quadrilaterals* was presented. These are quadrilaterals which have all four sidelengths integers, integral area, and integral diagonal lengths.

For $\delta = K8mn(m^2 - n^2)$, the family of trapezoids $O_2M_2M_1O_1$ obtained, are never Heron Quadrilaterals. They have all four sidelengths $|\overline{OO_1}|, x, y$, and $\frac{\alpha}{2}$ being integers, integral area, but the two diagonal lengths $d_1$ and $d_2$ are always irrational numbers.

Finally, for $\delta = K8mn(m^2 - n^2)$, the formulas in (11) and (12), and including $\alpha/2, \beta$, and $\gamma$; yield

$R_1 = K(m^2 - n^2)(m^2 + n^2)^2$, $R_2 = K \cdot 2mn(m^2 + n^2)^2$, $|\overline{O_1O_2}| = K \cdot (m^2 + n^2)^3$

Area of triangle $OO_1O_2 = K^2 \cdot mn \cdot (m^2 - n^2)(m^2 + n^2)^4$

$x = K(m^2 - n^2)^2 \cdot (m^2 + n^2)$ , $y = K \cdot 4 \cdot (mn)^2 \cdot (m^2 + n^2)$

Area of trapezoid $O_1O_2M_2M_1 = K^2 \cdot 2mn(m^2 - n^2)(m^2 + n^2)^4$

$$d_1 = K(m^2 - n^2)(m^2 + n^2)\sqrt{m^4 + 14m^2n^2 + n^4}$$

$$d_2 = K \cdot 4mn(m^2 + n^2)\sqrt{m^4 - m^2n^2 + n^4}$$

$\alpha/2 = K \cdot 4mn \cdot (m^2 - n^2)(m^2 + n^2)$, $\beta = K \cdot 16(mn)^2(m^2 - n^2)$

$\gamma = K \cdot 8mn(m^2 - n^2)^2$

## 7. Numerical Examples

**Table 1** (For triangle $AB\Gamma$)

| | $\alpha$ | $\beta$ | $\gamma$ |
|---|---|---|---|
| $K=1,\ m=2,\ n=1$ | 240 | 192 | 144 |
| $K=1,\ m=3,\ n=2$ | 3120 | 2880 | 1200 |
| $K=1,\ m=4,\ n=1$ | 8160 | 3840 | 7200 |

**Table 2**

| | $R_1$ | $R_2$ | $\lvert\overline{O_1O_2}\rvert$ | Area $OO_1O_2$ | $x$ | $y$ | $\alpha/2$ | $d_1$ | $d_2$ | Area $O_1O_2M_2M_1$ |
|---|---|---|---|---|---|---|---|---|---|---|
| $K=1, m=2, n=1$ | 75 | 100 | 125 | 3750 | 45 | 80 | 120 | $15\sqrt{61}$ | $40\sqrt{13}$ | 7500 |
| $K=1, m=3, n=2$ | 845 | 2028 | 2197 | 856830 | 325 | 1872 | 1560 | $65\sqrt{601}$ | $312\sqrt{601}$ | 1713660 |
| $K=1, m=4, n=1$ | 4335 | 2312 | 4913 | 5011260 | 3825 | 1088 | 4080 | $255\sqrt{481}$ | $272\sqrt{481}$ | 10022520 |

## 8. Closing Remarks

**Remark 1:** It easily follows from (3) and $\alpha^2 = \beta^2 + \gamma^2$; that

$$\left(\frac{1}{R_1}\right)^2 + \left(\frac{1}{R_2}\right)^2 = \left(\frac{4}{\alpha}\right)^2 = \left(\frac{1}{\alpha/4}\right)^2 . \quad (14)$$

Equation (14) implies the existence of a right triangle with leglengths $1/R_1$ and $1/R_2$; and hypotenuse length $\frac{1}{\alpha/4}$. Note that each of the line segments $\overline{\Gamma M_2}, \overline{M_2 O}, \overline{OM_1}$, and $\overline{M_1 B}$, in Figure 1, has length $\alpha/4$. Now, given a line segment $\overline{DE}$ of length $\ell$, one can construct, in standard Euclidean fashion, a line segment having length $1/\ell$, as illustrated in Figure 3: Just draw the line perpendicular to $\overline{DE}$ (standard Euclidean construction) and choose on it a point *F*, so that $\left|\overline{FE}\right| = 1$. Then draw the perpendicular at *F* to the line segment $\overline{DF}$, and let it intersect the straight line on which $\overline{DE}$ lies, at the point G. The three right triangles $EDF$, $EFG$, and $FDG$, are all similar; and from the similarity, it follows that

$$\frac{\left|\overline{FE}\right|}{\left|\overline{DE}\right|} = \frac{\left|\overline{EG}\right|}{\left|\overline{FE}\right|}; \quad \frac{1}{\ell} = \frac{\left|\overline{EG}\right|}{1} \quad \Leftrightarrow \quad \left|\overline{EG}\right| = 1/\ell .$$

Thus, given three line segments with lengths $R_{1,}$ $R_2$, and $\alpha$ mentioned above; we can construct three line segments having lengths; $\frac{1}{R_1}$, $1/R_2$, $\frac{1}{\alpha/4}$, and from there the above-mentioned right triangle.

One wonders what properties this right triangle may have. Note that it would be similar to $AB\Gamma$; since $\frac{R_1}{R_2} = \frac{\gamma}{\beta}$, by (3).

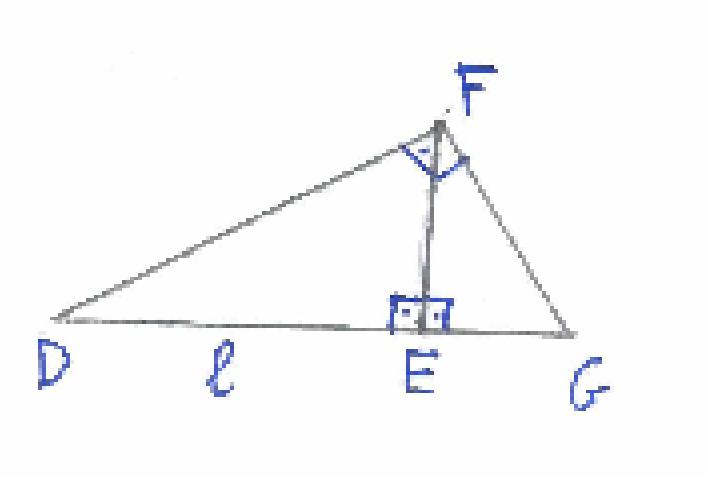

**Figure 3**

**Remark 2**: From Figure 1, $0 < \gamma < \beta \Leftrightarrow 0^\circ < \varphi < 45^\circ < \omega < 90^\circ$; and $\omega + \varphi = 90^\circ$. Also (from Figure 1 or by (3)) note that $0 < R_1 < \beta/2$ and $0 < \frac{\gamma}{2} < R_2$. By (3), it is clear that $R_1 = \gamma$ if, and only if, $R_2 = \beta$. A brief calculation shows that this occurs precisely when $\beta^2 + \gamma^2 = 4\beta\gamma \Leftrightarrow (\beta/\gamma)^2 - 4(\beta/\gamma) + 1 = 0$. Solving this quadratic equation yields, $\tan\omega = \beta/\gamma = 2 - \sqrt{3}$ or $2 + \sqrt{3}$. But $\tan\omega > 1$; and so, $\tan\omega = 2 + \sqrt{3} = \tan 75^\circ$ (see note below); $\omega = 75^\circ$ in this case.

The following can be readily established:

1) $0^\circ < \varphi < 45^\circ < \omega < 60^\circ \Leftrightarrow 0 < R_1 < R_2 < \gamma < \beta$.

2) $0^\circ < \varphi < 45^\circ < \omega = 60^\circ$; then $\varphi = 30^\circ$, $\gamma = \alpha/2$, $\beta = \frac{\sqrt{3}}{2}\alpha$; $0 < R_1 < R_2 = \gamma < \beta$. In this case, all four points $B, O_1, T$, and $O_2$; are aligned.

3) $0^\circ < \varphi < 45^\circ < 60^\circ < \omega < 75^\circ \Leftrightarrow 0 < R_1 < \gamma < R_2 < \beta$.

4) $0^\circ < \varphi < 45^\circ < 60^\circ < \omega = 75^\circ; \varphi = 15^\circ$. In this case $R_1 = \gamma$ and $R_2 = \beta$. The right triangles $AB\Gamma$ and $OO_1O_2$ are congruent.

5) $0^\circ < \varphi < 45^\circ < 60^\circ < 75^\circ < \omega < 90^\circ \Leftrightarrow 0 < \gamma < R_1 < \beta < R_2$.

*Note*: With regard to the $75^\circ$ angle. By using the double-angle identity for the tangent and $\tan 30^\circ = \frac{1}{\sqrt{3}}$; one can establish $\tan 15^\circ = 2 - \sqrt{3}$. Thus,

$$\tan 75^\circ = \cot 15^\circ = \frac{1}{\tan 15^\circ} = \frac{1}{2-\sqrt{3}} = \frac{2+\sqrt{3}}{\left(2-\sqrt{3}\right)\left(2+\sqrt{3}\right)} = \frac{2+\sqrt{3}}{4-3} = 2 + \sqrt{3}.$$

**References**

1. Rosen, Kenneth H., *Elementary Number Theory and its Applications*, 1993, Addison-Wesley Publishing Company (there is now a fourth edition as well), 544 p.p.
   ISBN 0-201-57889-1
   a) Four Pythagorean triples, see pages 436 – 442
   b) The result which states that the *n*th root of a positive integer is either an integer or otherwise irrational; follows from Th. 2.11 in that book, found on page 96.
   c) The result which states that if $a$ divides the product $bc$ and $(a,b)=1$; then $a$ must divide $c$; is stated on page 91 as Lemma 2.3.

2. Dickson, L. E., *History of Theory of Numbers, Vol. II*, AMS Chelsea Publishing, Rhode Island, 1992. ISBN: 0-8218-1935-6; 803 p.p. (unaltered text reprint of the original book, first published by Carnegie Institute of Washington in 1919, 1920, and 1923)
   a) For material on Pythagorean triangles and rational right triangles, see pages 165 – 190.
   b) For the Diophantine equation $x^4 + 14x^2y^2 + y^4 = z^2$, see page 635.

3. Sierpinski, W., *Elementary Theory of Numbers*, Warsaw, Poland, 1964, 480 p.p. (no ISBN number). More recent version (1988) published by Elsevier Publishing, and distributed by North-Holland. North-Holland Mathematical Library **32**, Amsterdam (1988). This book is available by various libraries, but it is only printed upon demand. Specifically, UMI Books on Demand from: Pro Quest Company, 300 North Zeeb Road, Ann Arbor, Michigan, 48106-1356 USA; ISBN: 0-598-52758-3

   a) For a description and derivation of Pythagorean triples, see pages 38 – 42.
   b) For the result which states that if an integer $a$ is a divisor $bc$ and $(a,b)=1$; then $a$ must be a divisor of $c$; see Theorem 5 on page 14.
   c) For the result which states that a positive integer is equal to the *n*th power of a rational number if, and only if, that positive integer is the *n*th power of an integer (which is equivalent to the statement that the *n*th root of a positive integer is either an integer or an irrational number), see Th. 7 on page 16.
   d) For the Diophantine equation $x^4 + x^2y^2 + y^4 = z^2$, see pages 73 – 74.

4. Pocklington, H. C., *Some Diophantine impossibilities*, Proc. Cambridge Philosophical Society, 17 (1914), p.p. 108 – 121

5. K. R. S. Sastry, *A description of a family of Heron Quadrilaterals*, Mathematics and Computer Education, Winter 2005, p.p. 72 -77